\input amstex
\input Amstex-document.sty

\def\b1{\text{\bf 1}}

\def\BZ{{\Bbb Z}}
\def\CA{{\Cal A}}
\def\CB{{\Cal B}}
\def\cC{\buildrel\vee\over{C}}

\def\CG{{\Cal G}}

\def\CO{{\Cal O}}

\def\CS{{\Cal S}}
\def\CT{{\Cal T}}
\def\CV{{\Cal V}}

\def\dpar{\partial}

\def\fA{{\frak A}}
\def\fAlg{{\frak{Alg}}}
\def\fb{{\frak b}}

\def\fU{{\frak U}}

\def\btu{\bigtriangleup}
\def\dplus{\buildrel{.}\over{+}}
\def\hra{\hookrightarrow}
\def\iso{\buildrel\sim\over\longrightarrow}

\def\lera{\leftrightarrow}

\def\lra{\longrightarrow}



\pageno 525

\def\shorttitle{Sur les Alg\`ebres Vertex Attach\'ees aux
Vari\'et\'es Alg\'ebriques}
\def\shortauthor{Vadim Schechtman}

\topmatter
\title\nofrills{\boldHuge Sur les Alg\`ebres Vertex Attach\'ees aux
Vari\'et\'es Alg\'ebriques}
\endtitle

\author \Large Vadim Schechtman* \endauthor

\thanks *Laboratoire Emile Picard, Universit\'e  Paul
Sabatier, 118, Route de Narbonne, 31062 Toulouse Cedex 4, France. E-mail: schechtman\@picard.ups-tlse.fr
\endthanks

\abstract\nofrills \centerline{\boldnormal Abstract}

\vskip 4.5mm

{\ninepoint {\it Sganarelle:}\ ... Mais encore faut-il croire quelque chose  dans le monde: qu'est-ce donc que
vous croyez?

{\it Dom Juan:}\ Ce que je crois?

{\it Sganarelle:}\ Oui.

{\it Dom Juan:}\ Je crois que deux et deux sont quatre, Sganarelle, et que
quatre et quatre sont huit.

{\it Moli\`ere, Dom Juan}

One dicusses sheaves of vertex algebras over smooth varieties and their connections with characteristic classes.

\vskip 4.5mm


}
\endabstract
\endtopmatter

\document

\baselineskip 4.5mm \parindent 8mm

\specialhead \noindent \boldLARGE 1. Introduction \endspecialhead

Le but de cette note est de presenter une classification de
certaines  alg\`ebres vertex, qui peuvent  \^etre  associ\'ees \`a
des vari\'et\'es alg\'ebriques lisses;  ceci est l'occasion de
rencontrer des classes caract\'eristiques ``style
Pontryagin-Atiyah-Chern-Simons". Ceci a \'et\'e  obtenu dans [GMS]
dont la presente note est un compl\'ement. On propose ici  une
d\'efinition plus  simple d'une {\it alg\'ebro\"{\i}de vertex}
({\it infra},  2.4, 2.6), un \'enonc\'e plus pr\'ecis et une
d\'emonstration courte du r\'esultat principal de {\it op. cit.}
({\it infra}, 3.4, 3.6, 3.7). \`A la fin on propose une
construction  directe des alg\`ebres vertex associ\'ees aux
courbes, \`a l'aide des alg\`ebres de Virasoro introduites par
A.Beilinson. Le point de d\'epart de cette note \'etait une
tentative de comprendre {\it le complexe de de~Rham chiral}
d\'ecouvert par F.Malikov, [MSV]. Je remercie vivement mes amis et
collaborateurs Arkady Vaintrob, Fyodor Malikov et Vassily
Gorbounov.

\specialhead \noindent \boldLARGE 2. Alg\'ebro\"{\i}des vertex
\endspecialhead

{\bf 2.1.} On fixe un corps de base $k$ de caract\'eristique $0$.
Rappelons qu'une {\it alg\`ebre vertex} est un $k$-espace
vectoriel $V$ muni d'un vecteur distingu\'e $\b1\in V$ (dit {\it
vacuum}) et d'une famille d'applications $k$-lin\'eaires $_{(i)}:\
V\otimes_kV\lra V\ (i\in\BZ)$ telles que $x_{(i)}y=0$ pour $i$
assez grand. Si l'on pose $\dpar x:=x_{(-2)}\b1$, on obtient un
op\'erateur $\dpar:\ V\lra V$. Les axiomes de [B] doivent \^etre
v\'erifi\'es. On n'est interess\'e que par les alg\`ebres vertex
{\it $\BZ_{\geq 0}$- gradu\'ees}, ce qui signifie que l'espace $V$
est muni d'une $\BZ_{\geq 0}$-graduation  (dite {\it poids
conforme}),  $V=\oplus_{n\geq 0}\ V_n$, $\b1\in V_0$ et
$V_{n(i)}V_m\subset V_{n+m-i-1}$. En particulier $\dpar V_n\subset
V_{n+1}$. Les morphismes de telles alg\`ebres \'etant d\'efinis de
mani\`ere \'evidente,  on obtient la cat\'egorie $\CV ert$ des
alg\`ebres vertex $\BZ_{\geq 0}$-gradu\'ees.

{\bf 2.2. ``Donn\'ees classiques" associ\'ees \`a une alg\`ebre
vertex.} Soit $V\in\CV ert$. Pour \^etre bref on \'ecrira $xy$ au
lieu de $x_{(-1)}y$; c'est une op\'eration non commutative et non
associative en g\'en\'eral. On a $V_nV_m\subset V_{n+m}$.   Posons
$A(V)=V_0$; l'op\'eration $xy$ est commutative et associative sur
$V_0$; donc $A(V)$ devient une $k$-alg\`ebre commutative avec
unit\'e $\b1$. Posons $\CA(V)=V_1$. Soit $\Omega(V)$ le
sous-$k$-vectoriel de $\CA(V)$ engendr\'e par les \'el\'ements
$a\dpar b,\ a,b\in A(V)$. Alors $\Omega(V)$ devient un
$A(V)$-module et $\dpar:\ A(V)\lra\Omega(V)$ est une d\'erivation.
En outre, si l'on pose $T(V):=\CA(V)/\Omega(V)$, l'op\'eration
$ax$ induit une structure de $A(V)$-module sur $T(V)$. (Par
contre, $\CA(V)$ n'est pas un $A(V)$-module en g\'en\'eral, \`a
cause de la non associativit\'e de l'op\'eration  $ax$.)
L'op\'eration $_{(0)}:\ \CA(V)\times\CA(V)\lra\CA(V)$ induit
l'application $[ , ]:\ T(A)\times T(A)\lra T(A)$ qui est un
crochet de Lie; l'op\'eration $_{(0)}:\ \CA(V)\times A(V)\lra
A(V)$ induit  une action de $T(V)$ sur $A(V)$ par d\'erivations;
on a  $[\tau, a\tau']=a[\tau,\tau']+\tau(a)[\tau,\tau']$, i.e.
$T(V)$ devient  une {\it $A(V)$-alg\'ebro\"{\i}de de Lie}. La
premi\`ere op\'eration induit aussi une action de l'alg\`ebre de
Lie $T(V)$ sur $\Omega(V)$ telle que $\dpar$ est un morphisme de
$T(V)$-modules, et   $\tau(a\omega)=\tau(a)\omega+a\tau(\omega)$.
Enfin, l'op\'eration $_{(1)}:\ \CA(V)\times\CA(V)\lra A$ est
sym\'etrique et induit un accouplement $A(V)$-bilin\'eaire
$\langle, \rangle:\ T(V)\times\Omega(V)\lra A(V)$ telle que
$\tau(\langle\tau',\omega\rangle)=\langle
[\tau,\tau'],\omega\rangle+ \langle\tau',\tau(\omega)\rangle$ et
$(a\tau)(\omega)=a\tau(\omega)+\langle\tau,\omega\rangle\dpar a$.

{\bf 2.3. ``Donn\'ees quantiques."} Les propri\'et\'es (Alg1) ---
(Alg3) ci-dessous sont v\'erifi\'ees, o\`u $a\in A(V),\
x,y,z\in\CA(V)$, $\pi:\ A(V)\lra T(A)$ \'etant la projection
canonique.

{\bf (Alg1)} $(ax)_{(1)}y=a(x_{(1)}y)-\pi(x)\pi(y)(a).$

{\bf (Alg2)} $x_{(0)}y + y_{(0)}x=\dpar(x_{(1)}y);\ (\dpar
x)_{(0)}y=0.$

{\bf (Alg3)}
$x_{(0)}(y_{(i)}z)=(x_{(0)}y)_{(i)}z+y_{(i)}(x_{(0)}z),\ i=0,1.$

{\bf 2.4.} Soit $A$ une $k$-alg\`ebre commutative de type fini, lisse sur $k$.
Posons $\Omega(A)=\Omega^1_{A/k},\ T(A)=Der_k(A,A)$ (l'alg\`ebre de Lie
de $k$-d\'erivations de $A$). Soit $\dpar=\dpar_{DR}:\ A\lra\Omega(A)$ la
d\'erivation  universelle. On a l'accouplement non degen\'er\'e
$A$-bilin\'eaire  $\langle,\rangle:\ T(A)\times\Omega(A)\lra A$; l'alg\`ebre
de Lie $T(A)$  agit sur $\Omega(A)$ par la d\'eriv\'ee de Lie. Ces donn\'ees
v\'erifient toutes  les propri\'et\'es de 2.2.

{\it Une $A$-alg\'ebro\"{\i}de vertex} est un $k$-espace vectoriel $\CA$
muni d'un sous-espace $F^1\CA\subset\CA$ avec les identifications de
$k$-vectoriels $F^1\CA=\Omega(A),\ \CA/F^1\CA=T(\CA)$ et des
op\'erations $k$-bilin\'eaires $_{(-1)}:\ A\times\CA\lra\CA,\ (a,x)\mapsto
ax=a_{(-1)}x$,  $_{(1)}:\ \CA\times\CA\lra A$ sym\'etrique, $_{(0)}:\
\CA\times\CA\lra\CA$.  On demande que (i) $A_{(-1)}\Omega(A)\subset \Omega(A)$
et que l'action de $A$  sur $\Omega(A)$ et sur $T(A)$ induite par
$_{(-1)}$
co\"{\i}ncide avec l'action canonique;  (ii) $\Omega(A)_{(i)}\Omega(A)=0\
(i=0,1)$; $\Omega_{(0)}\CA\subset\Omega(A)$, l'op\'eration $T(A)\times
T(A)\lra T(A)$  induite par $_{(0)}$ co\"{\i}ncide avec le crochet de Lie,
et l'action induite $T(A)\times\Omega(A)\lra\Omega(A)$ co\"{\i}ncide avec
la deriv\'ee de Lie; (iii) l'acccouplement $\langle,\rangle:\
T(A)\times\Omega(A)\lra A$ induit par  $_{(1)}$ co\"{\i}ncide avec
l'accouplement canonique.  Enfin, les propri\'et\'es (Alg1) --- (Alg3) doivent
\^etre v\'erifi\'ees.  Dans (Alg3) pour $i=1$ on
interpr\`ete la partie de gauche comme $\pi(x)(y_{(1)}z)$.

Soit $T^\cdot(A)$ une alg\'ebre de Lie dg concentr\'ee en degr\'es $-1, 0$,
avec $T^{-1}(A)=T^0(A)=T(A)$, $d:\ T^{-1}(A)\lra T^0(A)$ l'identit\'e, le
crochet $[,]_{0,-1}$ l'action adjointe. Soit $\Omega^\cdot(A):\
0\lra A\buildrel\dpar\over\lra\Omega(A)\lra\Omega(A)/\dpar A\lra 0$ le
complexe concentr\'e en degr\'es $-2, -1, 0$ avec les diff\'erentielles
\'evidentes. Ce complexe est un module dg sur $T^\cdot(A)$ (l'action de
$T^0(A)$ est  par la d\'eriv\'ee de Lie, la composante $[,]_{-1,-1}:\
T^{-1}(A)\times \Omega^{-1}(A)\lra\Omega^{-2}(A)$ \'etant l'accouplement
canonique, et  la composante $[ , ]_{-1,0}$ \'etant d\'efinie par
$[\tau,\bar{\omega}]= i_{\tau}(d\omega)$, o\`u $\bar{\omega}\in\Omega(A)/\dpar
A$ est l'image de  $\omega\in\Omega(A)$, $d:\
\Omega^1_{A/k}\lra\Omega^2_{A/k}$ est la diff\'erentielle de de~Rham,
$i_{\tau}:\ \Omega^2_{A/k}\lra\Omega^1_{A/k}$  est la convolution avec
$\tau$). On peut exprimer les axiomes (Alg2) et  (Alg3) en disant que l'on a
une alg\`ebre de Lie dg $\CA^\cdot:\  0\lra A\lra\ \CA\lra\CA/\dpar A\lra 0$,
concentr\'ee en degr\'es $-2, -1, 0$,  extension de $T^\cdot(A)$ par
$\Omega^\cdot(A)$ (consid\'er\'ee comme une  sous-alg\`ebre de Lie
ab\'elienne), telle que l'action de $T^\cdot(A)$ sur  $\Omega^\cdot(A)$ induite
co\"{\i}ncide avec celle d\'ecrite ci-dessus.
Un morphisme $g:\
\CA\lra\CA'$ est une application $k$-lin\'eaire respectant
les op\'erations $_{(i)}$ et les filtrations, qui
induit l'identit\'e sur $\Omega(A), T(A)$. D'o\`u la cat\'egorie $\CA lg_A$  des
$A$-alg\'ebro\"{\i}des vertex,  qui est {\it un  groupo\"{\i}de} (chaque
morphisme est un isomorphisme).

{\bf 2.5.} Soit $A$ comme dans 2.4. On d\'efinit la cat\'egorie $\CV ert_A$ dont
les objets sont $V\in\CV ert$ munies d'un isomorphisme de $k$-alg\`ebres
$A(V)\iso A$,  cet isomorphisme identifiant les donn\'ees classiques
$(T(V),\Omega(V),\dpar,  \langle,\rangle)$ correspondantes avec les donn\'ees
standardes  $(T(A),\Omega(A),\dpar_{DR},\langle,\rangle)$ d\'ecrites dans 2.4.
Les morphismes sont les morphismes des alg\`ebres vertex induisants
l'identit\'e sur les donn\'ees classiques.

La construction 2.2, 2.3 donne lieu au foncteur $Alg:\ \CV ert_A\lra
\CA lg_A,\ V\mapsto\CA(V)$.   Ce foncteur admet l'adjoint \`a gauche  $U:\ \CA
lg_A\lra\CV ert_A$, l'alg\`ebre vertex  $U\CA$ \'etant appel\'ee {\it
l'alg\`ebre enveloppante} d'un alg\'ebro\"{\i}de  vertex $\CA$. Pour chaque
$\CA\in\CA lg_A$ le morphisme d'adjonction  $\CA\lra Alg(U\CA)$ est un
isomorphisme.

{\bf 2.6.} Le langage suivant est un peu plus explicite et est
parfois commode.  Appelons {\it $A$-alg\'ebro\"{\i}de vertex
scind\'ee} un couple $\CB=\bigl(\langle,\rangle,c)$, o\`u
$\langle,\rangle:\ T(A)\times T(A)\lra A$ (resp. $c:\ T(A)\times
T(A)\lra\Omega(A)$) est une application $k$-bilin\'eaire
sym\'etrique (resp. antisym\'etrique). On demande que les
propri\'et\'es (AlgScind1)--(AlgScind3) ci-dessous soient
v\'erifi\'ees.

{\bf (AlgScind1)} $\langle a\tau,b\tau'\rangle -
a\langle\tau,b\tau'\rangle - b\langle a\tau,\tau'\rangle +
ab\langle\tau,\tau'\rangle=-\tau'(a)\tau(b).$

{\bf (AlgScind2)} $\langle\tau'',c(\tau,\tau')\rangle+ \langle\tau',c(\tau,\tau'')\rangle=\langle
[\tau,\tau'],\tau''\rangle + \langle\tau',[\tau,\tau'']\rangle - \tau(\langle\tau',\tau''\rangle)$
$+\tau'(\langle\tau,\tau''\rangle)/2 + \tau''(\langle\tau,\tau'\rangle)/2.$

{\bf (AlgScind3)} $3\bigl\{\tau(c(\tau',\tau''))+\tau'(c(\tau'',\tau))+\tau''(c(\tau,\tau'))-
c([\tau,\tau'],\tau'')-$ \linebreak $c([\tau',\tau''],\tau)-c([\tau'',\tau],\tau')\bigr\} =\dpar
\bigl\{\langle\tau,\frac{1}{2}[\tau',\tau'']+c(\tau',\tau'')\rangle
+\langle\tau',\frac{1}{2}[\tau'',\tau]+c(\tau'',\tau)\rangle
+\langle\tau'',\frac{1}{2}[\tau,\tau']+c(\tau,\tau')\rangle\bigr\}.$

Ses propri\'et\'es
entra\^{\i}nent que $\langle,\rangle$ et $c$ sont des op\'erateurs
diff\'erentiels, d'ordres $2$ et $3$ respectivement.

\'Etant donn\'es deux $A$-alg\'ebro\"{\i}des vertex scind\'ees
$\CB=(\langle,\rangle,c)$  et $\CB'=(\langle,\rangle',c')$, {\it
un morphisme} $f:\ \CB\lra\CB'$ est par d\'efinition une
application $k$-lin\'eaire $h=h_f:\ T(A)\lra\Omega(A)$
satisfaisant les propri\'et\'es (Mor1)--(Mor3) ci-dessous (dont la
premi\`ere implique que $h$ est un op\'erateur diff\'erentiel
d'ordre $2$).

{\bf (Mor1)} $\langle\tau',h(a\tau)\rangle -\langle
a\tau,\tau'\rangle + \langle a\tau,\tau'\rangle'=
a\bigl\{\langle\tau',h(\tau)\rangle - \langle\tau,\tau'\rangle +
\langle\tau',\tau\rangle'\bigr\}.$

{\bf (Mor2)} $\langle\tau,\tau'\rangle-\langle\tau,\tau'\rangle'=
\langle\tau,h(\tau')\rangle+\langle\tau',h(\tau)\rangle.$

{\bf (Mor3)}

$c(\tau,\tau')-c'(\tau,\tau')=\tau'(h(\tau))-\tau(h(\tau'))+
h([\tau,\tau'])+\dpar\{\langle\tau,h(\tau')\rangle-
\langle\tau',h(\tau)\rangle\}/2.$

La composition est d\'efinie par $h_{ff'}=h_f+h_{f'}$; l'identit\'e est
$h_{id}=0$. D'o\`u on obtient le groupo\"{\i}de $\CA lg\CS cind_A$ des
$A$-alg\'ebro\"{\i}des vertex scind\'ees. \'Etant donn\'e $\CB$ comme
ci-dessus, on pose $\CA(\CB)=T(A)\oplus\Omega(A)$ et d\'efinit les
op\'erations $_{(i)},\ i=-1,0,1$ par les formules
$a_{(-1)}\tau=(a\tau,-\gamma(a,\tau))$, o\`u $\gamma(a,\tau)\in\Omega(A)$
est d\'efini par
$$
\langle\tau',\gamma(a,\tau)\rangle=\langle a\tau,\tau'\rangle -
a\langle\tau',\tau\rangle + \tau\tau'(a)
$$
(l'axiome (AlgScind1) signifie que cette expression est $A$-lin\'eaire
en $\tau'$); $\tau_{(0)}\tau'=([\tau,\tau'],-c(\tau,\tau')+\frac{1}{2}
\dpar\langle\tau,\tau'\rangle$, $\tau_{(1)}\tau'=\langle\tau,\tau'\rangle$.
Ceci d\'efinit $\CA(\CB)\in\CA lg_A$. Si $f:\ \CB\lra\CB'$ est comme
ci-dessus, on d\'efinit le morphisme $g(f):\ \CA(\CB)\lra\CA(\CB')$ par
$g(f)(\tau)=(\tau,h_f(\tau))$. Ceci d\'efinit un foncteur
$\CA lg\CS cind_A\lra\CA lg_A$ qui est une \'equivalence des cat\'egories.

{\bf 2.7. Exemple.} Supposons que $A$ est telle que
$T(A)$ soit un $A$-module libre
et il existe une $A$-base $\fb=\{\tau_1,\ldots,\tau_n\}$ de $T(A)$ telle que
$[\tau_i,\tau_j]=0$ pour tous $i,j$. Nous appelerons telles alg\`ebres
{\it petites} est les bases $\fb$ {\it ab\'eliennes}.
On pose
$$
\langle
a\tau_i,b\tau_j\rangle_\fb=-b\tau_i\tau_j(a)-a\tau_j\tau_i(b)-
\tau_i(b)\tau_j(a), \eqno{(2.7)_{\langle,\rangle}}
$$
$$
c_\fb(a\tau_i,b\tau_j)=\frac{1}{2}\{\tau_i(b)\dpar\tau_j(a)-\tau_j(a)\dpar
\tau_i(b)\}+\frac{1}{2}\dpar\{b\tau_i\tau_j(a)-a\tau_j\tau_i(b)\}.
\eqno{(2.7)_c}
$$
Alors $\CB_\fb=(\langle,\rangle_\fb,c_\fb)$ est une
$A$-alg\'ebro\"{\i}de  vertex scind\'ee.

\specialhead \noindent \boldLARGE 3. Classification
\endspecialhead

{\bf 3.1.} Une alg\`ebre $A$ \'etant toujours comme dans 2.4, on
d\'efinit un groupo\"{\i}de  $\CG r(\Omega_A^{[2,3\rangle})$ dont
les objets sont les formes diff\'erentielles {\it ferm\'ees}
$\omega\in\Omega_{A/k}^{3,fer}$, avec
$$Hom_{\CG r(\Omega_A^{[2,3\rangle})}(\omega,\omega')=
\{\eta\in\Omega^2_{A/k}|d\eta=\omega-\omega'\}.$$
 La composition
des morphismes est l'addition de $2$-formes. L'addition des
$3$-formes induit une structure d'un {\it groupe ab\'elien en
cat\'egories} sur ce groupo\"{\i}de.

On remarque que si $\CA, \CA'$ sont deux $A$-alg\'ebro\"{\i}des vertex
avec le m\^eme espace sous-jacent la m\^eme op\'eration $_{(1)}$, alors
$x_{(0)}y-x_{(0)'}y\in\Omega(A)$; cet \'el\'ement ne d\'epend que des $\pi(x),
\pi(y)$, o\`u  $\pi:\ \CA\lra T(A)$ est l'application canonique, d'o\`u
l'application  $c_{\CA,\CA'}:\ T(A)\times T(A)\lra\Omega(A)$. De plus,
cette application est $A$-bilin\'eaire, et
$\omega_{\CA,\CA'}(\tau,\tau',\tau''):= \langle
\tau,c_{\CA,\CA'}(\tau',\tau'')\rangle$ est antisym\'etrique en $\tau, \tau',
\tau''$, donc $\omega_{\CA,\CA'}$ peut  \^etre consider\'ee comme une
$3$-forme diff\'erentielle, et cette forme  est ferm\'ee.

R\'eciproquement, \'etant donn\'e
$\CA=\CA lg_A$ et
$\omega\in\Omega^{3,fer}_{A/k}$, on d\'efinit
$\CA'=\CA\dplus\omega\in\CA lg_A$ ayant le m\^eme espace sous-jacent que
$\CA$ et la m\^eme op\'eration $_{(1)}$, avec $_{(0)'}=_{(0)}-\omega$.

Si $g:\ \CA\dplus\omega\lra\CA\dplus\omega'$ est un morphisme, alors
$(g - id)(\CA)\subset\Omega(A),\ (g - id)|_{\Omega(A)}=0$, donc $g - id$
induit une application $h_g:\ T(A)\lra\Omega(A)$. La fonction,
$\eta_g(\tau,\tau'):=\langle\tau,h_g(\tau)\rangle$ est antisym\'etrique
en $\tau, \tau'$ et $A$-bilin\'eaire, donc peut \^etre consider\'ee comme une
$2$-forme diff\'erentielle; on a $d\eta=\omega -\omega'$. Ceci induit une
bijection  $Hom_{\CA lg_A}(\CA\dplus\omega,\CA\dplus\omega')=
\{\eta\in \Omega^2_{A/k}|d\eta=\omega-\omega'\}$.
On a
$Hom_{\CA lg_A}(\CA,\CA')=Hom_{\CA lg_A}(\CA\dplus\omega,\CA'\dplus\omega)$.
Cela d\'efinit une Action
$$
\dplus:\ \CA lg_A\times\CG r(\Omega_A^{[2,3\rangle})\lra \CA lg_A.
\eqno{(3.1.1)}
$$

{\bf 3.2. Th\'eor\`eme.} {\it Si $A$ est petite} (voir 2.7), {\it alors le
groupo\"{\i}de $\CA lg_A$ est un Torseur  sous $\CG r(\Omega_A^{[2,3\rangle})$
par rapport \`a l'Action} (3.1.1),  {\it c'est \`a dire, pour chaque
$\CA\in\CA lg_A$ le  foncteur $\CA\dplus ?:\ \CG
r(\Omega_A^{[2,3\rangle})\lra\CA lg_A$ est  une \'equivalence.}

Par exemple, l'ensemble $\pi_0(\CA lg_A)$ des classes d'isomorphisme
de $\CA lg_A$ est un torseur sous $H^3_{DR}(A)$.
Grace \`a 2.6 le Torseur $\CA lg_A$ est non-vide pour $A$ petite.

{\bf 3.3.} Soient $A$ petite, et $\fb=\{\tau_i\},
\fb'=\{\tau'_i\}$ deux bases ab\'eliennes, d'o\`u les alg\'ebro\"{\i}des
scind\'ees $\CB_\fb, \CB_{\fb'}$; on a  $\tau'_i=\phi^{ij}\tau_j$
(la r\`egle de Einstein est
sous-entendue), $\phi=(\phi^{ij})\in GL_n(A)$ (pour \^etre bref, on
\'ecrit  $\fb'=\phi\fb$).
On d\'efinit une application  $h_{\fb',\fb}:\ T(A)\lra\Omega(A)$, comme
\'etant l'unique op\'erateur satisfaisant  (Mor1), tel que
$\langle\tau'_i,h_{\fb',\fb}(\tau'_j)\rangle=
-\frac{1}{2}\langle\tau'_i,\tau'_j\rangle_\fb$.
De plus, on d\'efinit une application $c_{\fb',\fb}:\ T(A)\times T(A)\lra
\Omega(A)$ comme \'etant l'unique op\'erateur tel que
$\CB_{\fb',\fb}:=(\langle,\rangle_\fb,c_{\fb',\fb})$ soit une
alg\'ebro\"{\i}de  vertex scind\'ee, et $h_{\fb',\fb}$ soit un morphisme
d'alg\'ebro\"{\i}des  scind\'ees $\CB_{\fb'}\lra\CB_{\fb',\fb}$. D'o\`u
la $3$-forme  $\alpha_{\fb',\fb}\in\Omega^{3,fer}_{A/k}$ telle que
$\CB_{\fb}=\CB_{\fb',\fb}\dplus\alpha_{\fb',\fb}$.
Si $\fb''=\{\tau''_i\}$ est la troisi\`eme base ab\'elienne, avec
$\tau''_i=\psi^{ij}\tau'_j$, on d\'efinit
la $2$-forme $\beta_{\fb'',\fb',\fb}:=h_{\fb'',\fb'}+h_{\fb',\fb}-
h_{\fb'',\fb}\in\Omega^2_{A/k}$.

{\bf 3.4. Th\'eor\`eme.} {\it $\beta_{\fb'',\fb',\fb}=\frac{1}{2}tr\{\phi^{-1}
\psi^{-1}d\psi d\phi\},\
\alpha_{\fb',\fb}=\frac{1}{6}tr\{(\phi^{-1}d\phi)^3\}$.}

{\it D\'emonstration.} Il resulte de (2.7) que
$$
c_\fb(\tau'_i,\tau'_j)=
\frac{1}{2}tr\bigl\{\phi^{-1}\tau'_j(\phi)\phi^{-1}\dpar\tau'_i(\phi)
- \phi^{-1}\tau'_i(\phi)\phi^{-1}\tau'_j(\phi)\phi^{-1}\dpar\phi -
(i\lera j)\bigr\}, \eqno{(3.4)_c}
$$
$$
\langle\tau'_i,\tau'_j\rangle_\fb=tr\bigl\{ -2\phi^{-1}\tau'_i\tau'_j(\phi) +
\phi^{-1}\tau'_i(\phi)\phi^{-1}\tau'_j(\phi)\bigr\}
\eqno{(3.4)_{\langle,\rangle}}
$$
d'o\`u, en utilisant (Mor1),
$$
h_{\fb',\fb}(\tau''_i)=tr\bigl\{\phi^{-1}\dpar\tau''_i(\phi)-\frac{1}{2}
\phi^{-1}\tau''_i(\phi)\phi^{-1}\dpar\phi +
\phi^{-1}\psi^{-1}\tau''_i(\psi)\dpar\phi\bigr\}. \eqno{(3.4)_h}
$$
Par d\'efinition,
$c_{\fb',\fb}(\tau'_i,\tau'_j)=\tau'_i(h(\tau'_j))-
\tau'_j(h(\tau'_i));\ \ \alpha_{\fb',\fb}=c_\fb - c_{\fb',\fb}$,
d'o\`u;
$$
\alpha_{\fb',\fb}(\tau'_i,\tau'_j)=-\frac{1}{2}tr\bigl\{
\phi^{-1}\tau'_i(\phi)\phi^{-1}\tau'_j(\phi)\phi^{-1}\dpar\phi -
(i\lera j)\bigr\}. \eqno{(3.4)_\alpha}
$$
En outre, $(3.4)_h$ entra\^{\i}ne
$$
\beta_{\fb'',\fb',\fb}(\tau''_i)=\frac{1}{2}tr\bigl\{
\phi^{-1}\psi^{-1}\tau''_i(\psi)\dpar\phi -
\tau''_i(\phi)\phi^{-1}\psi^{-1}\dpar\psi
\bigr\}
\eqno{(3.4)_\beta}
$$
d'o\`u le th\'eor\`eme. Ici l'on identifie une $3$-forme $\alpha$ avec
une application antisym\'etrique $T(A)\times T(A)\lra\Omega(A)$ d\'efinie
par $\alpha(\tau,\tau')=i_{\tau}i_{\tau'}\alpha$.
$\btu$

{\bf 3.5. Classe de Pontryagin.} Soit $X$ une vari\'et\'e alg\'ebrique lisse
sur $k$, $E$ un fibr\'e vectoriel sur $X$. Choisissons une recouvrement
affine $\fU=\{ U_i\}$ de $X$, et des bases
$\fb^i$ des $\Gamma(U_i,\CO_X)$-modules $\Gamma(U_i,E)$, d'o\`u
le cocycle de $\cC$ech $\phi=(\phi_{ij}),\
\phi_{ij}\in\Gamma(U_{ij},GL_n(\CO_X))$, $\fb^i=\phi_{ij}\fb^j$ sur $U_{ij}$,
$\phi_{ij}\phi_{jk}=\phi_{ik}$ sur $U_{ijk}$. Consid\'erons les
cocha\^{\i}nes de $\cC$ech
$p_2(\phi)=\bigl(\frac{1}{2}
tr\{\phi_{jk}^{-1}\phi_{ij}^{-1}d\phi_{ij}d\phi_{jk}\}\bigr)\in
C^2(\fU,\Omega^2_X)$,
$p_3(\phi)=\bigl(\frac{1}{6}tr\{(\phi_{ij}^{-1}d\phi_{ij})^3\}\bigr)\in
C^1(\fU,\Omega^3_X)$; on a $d_{\cC ech}p_2(\phi)=0$,
$d_{DR}p_2(\phi)=d_{\cC ech}p_3(\phi), d_{DR}p_3(\phi)=0$.
Il en r\'esulte que
$p(\phi):=\bigl(p_2(\phi),p_3(\phi)\bigr)\in
Z^2(\fU,\Omega^{[2,3\rangle}_X)$
o\`u $\Omega^{[2,3\rangle}_X:=
\bigl(\Omega^2_X\lra\Omega^{3,fer}_X\bigr)$, la diff\'erentielle totale
dans le bicomplexe de $\cC$ech \`a coefficients dans ce complexe
\'etant $d=d_{DR}+(-1)^{|DR|}d_{\cC ech}$.
De plus, si l'on choisit des autres bases $'\fb^i=g_i\fb^i$, d'\`u
$g=(g_i)\in C^0(\fU,GL_n(\CO_X))$, le cocycle correspondant est $\phi'=
^g\phi$, o\`u $^g\phi_{ij}=g_i\phi_{ij}g_j^{-1}$. On d\'efinit
$$
p_2(\phi,g):=\bigl(\frac{1}{2}tr\bigl\{\phi_{ij}^{-1}g_i^{-1}dg_i\phi_{ij}
g_j^{-1}dg_j+\phi_{ij}^{-1}d\phi_{ij}g_j^{-1}dg_j-g_i^{-1}dg_i
d\phi_{ij} \phi_{ij}^{-1}\bigr\}\bigr),
$$
$$
p_3(g):=\bigl(\frac{1}{6}tr\bigl\{(g_i^{-1}dg_i)^3\bigr\}\bigr);\
p(\phi,g)=\bigl(p_2(\phi,g),p_3(g)\bigr)\in
C^1(\fU,\Omega_X^{[2,3\rangle}).
$$
Alors $p_3(^g\phi)=p_3(\phi)+d_{\cC ech}p_3(g)+
d_{DR}p_2(\phi,g), p_2(\phi^g)=p_2(\phi)+d_{\cC ech}p_2(\phi,g)$, d'o\`u
$p(^g\phi)=p(\phi)+dp(\phi,g)$.
Donc la classe $p(E)$ de
$p(\phi)$ dans  $H^2(X,\Omega^{[2,3\rangle}_X)$ qu'on  peut appeler {\it la
classe de  Pontryagin-Atyiah-Chern-Simons (pacs)}, ne depend que de $E$. On
remarque que $p(\phi)=p(\phi^{-1t})$, donc  $p(E)=p(E^*)$.

{\bf 3.6.} Les groupo\"{\i}des $\CA lg_{\Gamma(U,\CO_X)},\ U\subset X,$
forment  un champ $\fAlg_X$ sur la topologie de Zariski (m\^eme \'etale),
parsque les op\'erations $_{(i)}$ sont des op\'erateurs diff\'erentiels qui
se localisent. D'apr\`es
3.2, $\fA lg_X$ est une {\it gerbe} sous $\Omega^{[2,3\rangle}_X$  (localement
non-vide, mais pas localement connexe). Donc {\it la classe caract\'eristique}
$c(\fAlg_X) \in H^2(X,\Omega_X^{[2,3\rangle})$ est d\'efinie, telle que
$c(\fAlg_X)=0$ ssi $\Gamma(X,\fA lg_X)$ est non-vide.  Rappelons sa
d\'efinition. On choisit un recouvrement affine $\fU=\{U_i\}$  de $X$ avec
$U_i$ petites;
on choisit les objets $\CA_i\in \Gamma(U_i,\fAlg_X)$. Sur les doubles
intersections, il existe les isomorphismes $h_{ij}:\ \CA_j|_{U_{ij}}\iso
\CA_i|_{U_{ij}} \dplus\alpha_{ij}$, $\alpha_{ij}\in\Omega^{3,fer}(U_{{ij}})$.
Si l'on pose
$\beta_{ijk}:=h_{ij}|_{U_{ijk}}-h_{ik}|_{U_{ijk}}+h_{jk}|_{U_{ijk}}\in
\Omega^2(U_{ijk})$, on a
$c(\{\CA_i\},\{h_{ij}\}):=(\bigl(\alpha_{ij}),(\beta_{ijk})\bigr)\in
Z^2(\fU,\Omega_X^{[2,3\rangle})$. Pour une autre famille $(\{\CA'_i\},
\{h'_{ij}\})$ il existent $h_i:\ \CA_i\iso\CA'_i\dplus\alpha_i,\
\alpha_i\in\Omega^{3,fer}_{A/k}$; alors $(h_j\dplus\alpha_{ij})\circ h_{ij}:\
\CA_i\iso\CA'_j\dplus(\alpha_j+\alpha_{ij})$ et
$(h'_{ij}\dplus\alpha_i)\circ h_i:\ \CA_i\iso\CA'_j\dplus (\alpha'_{ij}+
\alpha_i)$, donc il existe l'unique $\beta_{ij}\in\Omega^2_{A/k}$ telle que
$d\beta_{ij}=\alpha'_{ij}-\alpha_{ij}+\alpha_i-\alpha_j$. Alors
$d((\alpha_i),(\beta_{ij}))=c(\{\CA'_i\},\{h'_{ij}\})-c(\{\CA_i\},\{h_{ij}\})$;
par d\'efinition $c(\fAlg_X)$ est la classe de $c(\{\CA_i\},\{h_{ij}\})$
dans la cohomologie.

Soit $\CT_X$ le fibr\'e tangent de $X$.
Choisissons des bases bonnes $\fb^i$ de $\Gamma(U_i,\CT_X)$, avec
$\fb^i=\phi_{ij}\fb^j$, $\phi=(\phi_{ij})\in Z^1(\fU,GL_n(\CO_X))$.
Alors, d'apr\`es 3.4, $\alpha_{\fb^i\fb^j}=p_3(\phi)_{ij}$ et
$\beta_{\fb^i\fb^j\fb^l}:=h_{\fb^j\fb^l}-h_{\fb^i\fb^l}+h_{\fb^i\fb^j}=
p_2(\phi)_{ijl}$. De plus, si $\{'\fb^i\}$ est une autre famille des bases
bonnes, avec $'\fb^i=g_i\fb^i,\ g=(g_i)$, alors $\alpha_{'\fb^i\fb^i}=
p_3(g)_i$ et $h_{'\fb^i\
'\fb^j}-h_{\fb^i\fb^j}+h_{'\fb^i\fb^i}-h_{'\fb^j\fb^j}= p_2(\phi,g)_{ij}$. En
particulier, on a

{\bf 3.7. Th\'eor\`eme.} {\it $c(\fAlg_X)=p(\CT_X)$, o\`u $\CT_X$ est
le fibr\'e tangent de $X$.} $\btu$

Soit $\phi$ comme dans 3.6, $p=p(\phi)$; soit
$\CG r_p$ le
groupo\"{\i}de dont les objets sont les $1$-cocha\^{\i}nes de Cech
$\omega\in C^1(\fU,\Omega_X^{[2,3\rangle})$ telles que $d\omega=p$,
avec $Hom_{\CG r_p}(\omega,\omega')=\{\eta\in
C^0(\fU,\Omega_X^{[2,3\rangle})|\ d\eta=\omega-\omega'\}$. La construction
3.6 donne lieu au foncteur $\CG_p\lra \Gamma(X,\fAlg_X)$ qui
est une \'equivalence des cat\'egories. Il en r\'esulte que
$\pi_0\Gamma(X,\fAlg_X)$ est un torseur sous $H^1(X,\Omega_X^{[2,3\rangle})$,
non-vide si et seulement si $p(\CT_X)=0$, et pour $\CA\in\Gamma(X,\fAlg_X)$
le groupe $Aut (\CA)$ est isomorphe \`a $H^0(X,\Omega_X^{[2,3\rangle})$.

\specialhead \noindent \boldLARGE 4. Exemple \endspecialhead

Soit $X$ une courbe lisse sur $k$. Dans ce cas
$\Omega_X^{[2,3\rangle}=0$, donc sur $X$ il existe l'unique, \`a
isomorphisme unique pr\`es, $\CO_X$-alg\'ebro\"{\i}de vertex
$\CA_X$. On propose ici une construction directe de $\CA_X$. Pour
$j\in\BZ$ on a d\'efini dans [BS] le faisceau d'alg\`ebres de Lie
diff\'erentielles  gradu\'ees $\CA_j^{\cdot}$ ({\it $j$-i\`eme
Virasoro}) sur $X$ (cf. [BS] 3.1). On a $\CA_j^i=0$ pour $i\neq
-2,-1,0$; $\CA_j^{-2}=\CO_X,\
  \CA_j^0=\CT_X$. On a la suite exacte canonique des $k$-vectoriels (des
$\CO_X$-modules si $j=0$)
$0\lra\Omega^1_X\lra\CA_j^{-1} \buildrel{\pi}\over\lra\CT_X\lra 0$;
Par
d\'efinition, la diff\'erentielle $d: \CA_j^{-1}\lra\CA_j^0$ est  \'egale \`a
$\pi$ et $d:\ \CA_j^{-2}\lra \CA_j^{-1}$ est \'egale \`a la compos\'ee de la
diff\'erentielle  de de Rham avec l'inclusion $\Omega^1_X\hra\CA_j^{-1}$.
Comme il est expliqu\'e dans {\it op. cit.}, la cat\'egorie des alg\`ebres
de Lie dg comme ci-dessus est un {\it $k$-espace vectoriel en cat\'egories};
en particulier, on peut les multiplier par un scalaire. On a l'isomorphisme
canonique $\CA^\cdot_j\iso (6j^2-6j+1)\CA^\cdot_0$. Pour chaque
$\lambda\in k$ on a l'isomorphisme des $k$-modules
$\lambda\CA^{-1}_0\iso\CA^{-1}_0$, donc la structure canonique d'un
$\CO_X$-module sur $\CA_0^{-1}$ induit une structure de $\CO_X$-module sur
$\lambda\CA_0^{-1}$.

Consid\'erons l'alg\`ebre de Lie dg $6\CA_0$. On pose
$\CA_X=6\CA_0^{-1}$. On d\'efinit les op\'erations par
$a_{(-1)}x=ax-2\dpar\pi(x)(a);\ x_{(0)}y=[\pi(x),y],\
x_{(1)}y=[x,y]\ (a\in \CO_X,\ x,y\in\CA_X)$. Alors les axiomes
(Alg1)--(Alg3) sont v\'erifi\'es, et l'on obtient une
alg\'ebro\"{\i}de vertex sur $X$.

\specialhead \noindent \boldLARGE Bibliographie \endspecialhead

\widestnumber\key{[XXXX]}

\ref \key BS \by A.~Beilinson, V.~Schechtman \paper \rm
Determinant bundles and Virasoro algebras \jour {\it Comm. Math.
Phys.}, \vol \rm 118 \yr 1988 \pages 651--701
\endref

\ref \key B \by R.~Borcherds \paper \rm Vertex algebras, Kac-Moody
algebras, and the Monster \jour {\it Proc. Natl. Acad. Sci, USA},
\vol \rm 83 \yr May 1986 \pages 3068--3071
\endref

\ref \key GMS \by V.~Gorbounov, F.~Malikov \& V.~Schechtman \book
\it Gerbes of chiral differential operators. II \publ
math.AG/0003170
\endref

\ref \key MSV \by F.~Malikov, V.~Schechtman \& A.~Vaintrob \paper
\rm Chiral de Rham complex \jour {\it Comm. Math. Phys.}, \vol \rm
204 \yr 1999 \pages 439--473
\endref

\enddocument